\documentclass[amscd,amssymb,11pt]{amsart}

\usepackage[all]{xy}

\newtheorem{thm}{Theorem}[section]
\newtheorem{lem}[thm]{Lemma}
\newtheorem{cor}[thm]{Corollary}
\newtheorem{prop}[thm]{Proposition}

\theoremstyle{definition}

\newtheorem{ex}[thm]{Example}

\theoremstyle{remark}

\newtheorem{rem}[thm]{Remark}

\numberwithin{equation}{section}

\newcommand{\thmref}[1]{Theorem~\ref{#1}}
\newcommand{\corref}[1]{Corollary~\ref{#1}}
\newcommand{\secref}[1]{\S\ref{#1}}

\newcommand{\propref}[1]{Proposition~\ref{#1}}
\newcommand{\lemref}[1]{Lemma~\ref{#1}}

\newcommand{\exref}[1]{Example~\ref{#1}}

\newcommand{\hocolim}{\operatorname*{hocolim}}

\newcommand{\colim}{\operatorname*{colim}}

\newcommand{\Map}{\operatorname{Map}}
\newcommand{\MapS}{\operatorname{Map_{\mathcal S}}}
\newcommand{\MapT}{\operatorname{Map_{\mathcal T}}}

\newcommand{\Tor}{\operatorname{Tor}}

\newcommand{\Sp}{{\mathcal  S}}

\newcommand{\T}{{\mathcal  T}}

\newcommand{\Z}{{\mathbb  Z}}

\newcommand{\R}{{\mathbb  R}}
\newcommand{\F}{{\mathbb  F}}

\newcommand{\Sinfty}{\Sigma^{\infty}}

\newcommand{\sm}{\wedge}

\newcommand{\ra}{\rightarrow}
\newcommand{\xra}{\xrightarrow}

\newcommand{\hra}{\hookrightarrow}

\begin{document}

\title[Mapping spaces and homology isomorphisms]{Mapping spaces and homology isomorphisms}             

\author[N.J.Kuhn]{Nicholas J.~Kuhn}                                  
\address{Department of Mathematics \\ University of Virginia \\ Charlottesville, VA 2290}    
\email{njk4x@virginia.edu}
\thanks{This research was partially supported by a grant from the National Science Foundation}     
 
\date{July 8, 2004.}


\begin{abstract}   Let $\Map(K,X)$ denote the space of pointed continuous maps from a finite cell complex $K$ to a space $X$.  Let $E_*$ be a generalized homology theory.  We use Goodwillie calculus methods to prove that under suitable conditions on $K$ and $X$, $\Map(K, X)$ will send a $E_*$--isomorphism in either variable to a map that is monic in $E_*$ homology.  Interesting examples arise by letting $E_*$ be $K$--theory, $K$ be a sphere, and the map in the $X$ variable be an exotic unstable Adams map between Moore spaces.
\end{abstract}
                                                  
\maketitle

\section{Introduction and main results} \label{introduction}

Let $K$ and $X$ be pointed spaces, with $K$ homotopy equivalent to a finite cell complex, and then let $\Map(K,X)$ denote the space of pointed continuous maps from $K$ to $X$.  Fixing $K$, this includes many important constructions on $X$.  For example, $\Map(S^n,X) = \Omega^n X$, the $n^{th}$ loopspace of $X$, and $\Map(S^1_+,X) = \mathcal LX$, the free loopspace on $X$.

Suppose $E_*$ is a generalized homology theory.  A fundamental problem is to try to determine to what extent $E_*(\Map(K,X))$ might be determined by $E_*(X)$.  This is difficult, and has a long history, even when $E_*$ is ordinary homology with field coefficients and $K=S^1$.

We consider a related problem.  A map  $f: X \ra Y$ will be called an $E_*$--isomorphism if $E_*(f)$ is an isomorphism.  One can ask to what extent does $\Map(K, \text{\hspace{.1in}})$ preserve $E_*$--isomorphisms? 

This is question of interest when $E_*$ is a nonconnective theory as the following simple example illustrates: the constant map $c: K(\mathbb Z/p, 2) \ra *$ is an isomorphism in complex $K$--theory $K_*$, but $\Omega c: K(\mathbb Z/p,1) \ra *$ is not.  A much more subtle family of examples has been constructed recently by  L. Langsetmo and D. Stanley \cite{langsetmostanley}: see \exref{exotic maps} below.

In this paper we use Goodwillie calculus methods to prove the curious result that under suitable conditions on $K$ and $X$, $\Map(K, X)$ will send an $E_*$--isomorphism in either variable to a map that is monic in $E_*$ homology.  

\subsection{The main theorem}

We need to define some numerical invariants of spaces.

Let $d(K)$ be the minimal $d$ such that $K$ is homotopy equivalent to  a $d$--dimensional complex. 

Let $e(K)$ be the minimal $n$ such that there exists a parallelizable $n$--dimensional manifold $M$, together with a closed subcomplex $A$ such that $K$ is homotopy equivalent to $M/A$. For example, $e(S^n) = n$, as $S^n = D^n/S^{n-1}$. 

Let $c(X)$ be  the connectivity of $X$.

Let $s(X)$ be the minimal $n$ such that $X$ is homotopy equivalent to an $n$--fold suspension.  
 
Armed with these definitions, we can state our main theorem.

\begin{thm} \label{main theorem}  Suppose $e(K) \leq s(X)$ and $d(K) \leq c(X)$.  \\

\noindent{\bf (1)} \ If $f: X \ra Y$ is an $E_*$--isomorphism, then 
$$ \Map(K,f): \Map(K,X) \ra \Map(K,Y)$$
is $E_*$--monic. \\

\noindent{\bf (2)} \ If $g: L \ra K$ is an $E_*$--isomorphism of finite complexes, then 
$$ \Map(g,X): \Map(K,X) \ra \Map(L,X)$$
is $E_*$--monic. 

\end{thm}

\begin{cor} \label{main cor} If $Z$ is connected, and $f: \Sigma^n Z \ra Y$ is an $E_*$--isomorphism, then $\Omega^n f: \Omega^n \Sigma^n Z \ra \Omega^n Y$ is $E_*$--monic. 
\end{cor}

\begin{rem} This corollary seems to be new even when $n=1$.  To the best of the author's knowledge, the only results of this sort in the literature are the author's papers \cite{k1,k2} which contain the $n=\infty$ version of the corollary.
\end{rem}

Note that $d(K) \leq e(K)$. Furthermore $s(X) \leq c(X)+1$ is always true, and very often $s(X) \leq c(X)$.  For example, $s(M^n(d)) = c(M^n(d)) = n-2$ where $M^n(d)$ is the Moore space $D^n \cup_d S^{n-1}$. Thus when the first inequality in the hypotheses of the theorem holds, so usually does the second.   In general, $e(K)$ seems hard to compute exactly.  The appendix includes some observations of Greg Arone and the author which yield some further explicit calculations, and some general bounds.  For example, $e(M^n(d)) = n+1$, and, $e(K) \leq 2d(k)-1$ for all $K$ with $d(K) \geq 1$.

The numeric hypotheses of our theorem are easy to explain.  The condition $d(K) \leq c(X)$ guarantees the strong convergence of the Goodwillie tower of the functor sending a space $X$ to the suspension spectrum $\Sinfty \Map(K,X)$.  The condition $e(K) \leq s(X)$ implies that there is a filtered configuration space approximation to $\Map(K,X)$, as in work of B\"odigheimer \cite{bodigheimer}, following McDuff \cite{mcduff} and May \cite{may}.

When both numeric conditions hold, statement (1) of \thmref{main theorem} is proved by using properties of Goodwillie towers to play the two corresponding geometric conditions against each other.  Using the existence of Bousfield localization of spaces, statement (2) is then a formal consequence of (1).

The proof of \thmref{main theorem} is given in sections 2 and 3.  In section 2, we outline how general calculus theory leads to theorems like ours, while in the shorter section 3, we specialize to the case in hand. 

\subsection{Examples and applications}

\begin{ex} \label{exotic maps}
Let $p$ be an odd prime.  For each $(m,n)$ in an explicit infinite list of pairs, with $m \geq 4$ and both $m$ and $(n-m)$ taking on arbitrarily large values, Langsetmo and Stanley \cite{langsetmostanley} construct a $K_*$--isomorphism 
$$ f: M^{n}(p) \rightarrow M^m(p)$$
such that $\Omega f$ is {\em not} a $K_*$--isomorphism. For example, with $p=3$, for all $t \geq 1$, one has such a nondurable $K_*$--isomorphism  $f: M^{4t}(3) \ra M^4(3).$  

It is not hard to deduce that then, for all $j \geq 1$, the 3--connected cover of $\Omega^j f$, 
$$(\Omega^j f)\langle 3 \rangle: \Omega^jM^{n}(p)\langle 3 \rangle \ra \Omega^jM^{m}(p)\langle 3 \rangle, $$ is also not a $K_*$--isomorphism\footnote{This follows from  a theorem of Bousfield \cite[Theorem 11.10]{bousfield3}, but is easy to prove directly, using that $\tilde{K}_*(K(\Z/p, 2) = 0$}.

In contrast, \corref{main cor}  implies that for all $1 \leq j \leq n-2$, 
$$\Omega^j f: \Omega^j M^{n}(p) \ra \Omega^j M^m(p)$$ is $K_*$--monic.  

Combining these results, we conclude that for $1 \leq j \leq n-5$, 
$$(\Omega^jf) \langle 3 \rangle : \Omega^j M^{n}(p) \ra \Omega^j M^m(p)\langle 3 \rangle$$ is $K_*$--monic but not $K_*$--epic.
\end{ex}

\begin{ex}  When the homology theory $E_*$ is $K(r)_*$, the $r^{th}$ Morava $K$--theory at a prime $p$, the corollary has the following computational implication.  

Let $f: \Sigma^n Z \ra Y$ be a $K(r)_*$--isomorphism, with $Z$ connected and $n \geq 1$, and let $F$ be the fiber of $f$.  The $K(r)_*$ bar spectral sequence associated to the principal fibration 
$$ \Omega^n \Sigma^n Z \xra{\Omega^n f} \Omega^n Y \ra \Omega^{n-1} F $$
converges to $K(r)_*(\Omega^{n-1}F)$ and has 
$$E^2_{*,*} = \Tor^{K(r)_*(\Omega^n \Sigma^n Z)}_{*,*}(K(r)_*(\Omega^n Y), K(r)_*).$$
By the corollary, 
$$ K(r)_*(\Omega^n \Sigma^n Z) \xra{(\Omega^n f)_*} K(r)_*(\Omega^n Y)$$
is monic.  The map $(\Omega^n f)_*$ is in the category of $\mathcal K/p$--Hopf algebras studied by Bousfield in \cite[Appendix]{bousfield2}.  He shows \cite[Thm.10.8]{bousfield2} that objects in this category are flat over subobjects, when viewed as algebras.  We conclude that the spectral sequence collapses, giving an isomorphism
$$ K(r)_*(\Omega^{n-1}F) \simeq K(r)_*(\Omega^n Y) \otimes_{K(r)_*(\Omega^n \Sigma^n Z)} K(r)_*$$
of $K(r)_*$--coalgebras\footnote{$K(r)_*$--Hopf algebras if $n >1$.}.
\end{ex}

\begin{ex}  Suppose $g: L \ra K$ is a $K(r)_*$--isomorphism between finite complexes.  Let $C$ be the cofiber of $g$.  Applying statement (2) of \thmref{main theorem} to $\Sigma g$, and reasoning as in the last example, we deduce that, for all $X$ such that $e(K) < s(X)$ and $d(K) < c(X)$, one gets an isomorphism of $K(r)_*$--coalgebras
$$ K(r)_*(\Map(C, X)) \simeq K(r)_*(\Map(\Sigma L, X))\otimes_{K(r)_*(\Map(\Sigma K, X))} K(r)_*. $$
\end{ex}

\subsection{Acknowledgements}

The simple argument given in \secref{section 3} proving that statement (1) of \thmref{main theorem} implies statement (2) is due to Pete Bousfield, and replaces a different argument by the author, which needed the side hypothesis that $E_*$ be a ring theory.  For this, and for other encouraging `e'--conversations, I offer Pete my thanks.  Thanks are also due my colleagues Greg Arone and Slava Krushkal for discussions about this material.

\section{Goodwillie calculus and $E_*$--isomorphisms} \label{section 2}

Let $\T$ denote the category of based spaces, and $\Sp$ a nice model category of spectra, e.g. the category of $S$--modules of \cite{ekmm}.  In this section we find conditions on a functor $F: \T \ra \Sp$ and a space $X$ ensuring that if $f: X \ra Y$ is an $E_*$--isomorphism, then $F(f): F(X) \ra F(Y)$ will be $E_*$--monic.

\subsection{Review of Goodwillie calculus}

 In the series of papers \cite{goodwillie1, goodwillie2, goodwillie3}, Tom Goodwillie has developed his theory of polynomial resolutions of homotopy functors.  We need to summarize some aspects of Goodwillie's work as they apply to functors $F: \T \ra \Sp$.  

As carefully discussed in \cite{goodwillie2, goodwillie3}, a functor is said to be polynomial of degree $r$ if it takes strongly homotopy cocartesion $(r+1)$--cubical diagrams to homotopy cartesian cubical diagrams.  In \cite{goodwillie3}, given a functor $F$ from one topological model category to another, Goodwillie proves the existence of a tower $\{P_rF\}$ under $F$ so that $F \ra P_rF$ is the universal arrow to a polynomial functor of degree $r$, up to weak equivalence. 

The functors $F$ of interest to us in this paper satisfy an additional property: they will be {\em finitary}. Here, following \cite[Definition 5.10]{goodwillie3}, $F$ is said to be finitary if it commutes with filtered homotopy colimits up to equivalence.  

Examination of the construction of $P_rF$ shows that $P_r$ satisfies the following useful properties.

\begin{lem} \label{Pr lemma} (Compare with \cite[Proposition 1.7]{goodwillie3}.)

\noindent{\bf (1)} If $F(X) \ra G(X) \ra H(X)$ is a fibration sequence for all $X$, so is $P_rF(X) \ra P_rG(X) \ra P_rH(X)$. \\

\noindent{\bf (2)} Given natural transformations $F_1 \ra F_2 \ra \dots$, the natural map 
$$\hocolim_s P_rF_s(X) \ra P_r(\hocolim_s F_s)(X)$$
 is an equivalence for all $r$ and $X$. \\

\noindent{\bf (3)} If $F$ is finitary, so is $P_rF$ for all $r$.

\end{lem}

The fact that the suspension of a strongly homotopy cocartesian cube is again strongly homotopy cocartesion implies the next property of Goodwillie towers.

\begin{lem} \label{Pr susp lemma} \cite[Remark 1.1]{goodwillie3}
There is a natural equivalence $$P_r(F \circ \Sigma^d)(X) \simeq (P_rF)(\Sigma^d X).$$
\end{lem} 

Let $D_rF(X)$ be the homotopy  fiber of $P_rF(X) \ra P_{r-1}F(X)$.  $D_rF$ is homogeneous of degree $r$: it has degree $r$, and $P_{r-1}D_rF(X)$ is weakly contractible. (This follows from \lemref{Pr lemma}(1): see \cite[Proposition 1.17]{goodwillie3}.) Goodwillie analyzes $D_rF$.  We need his description when $F$ is also finitary and takes values in a stable model category like $\Sp$.

\begin{prop} \cite[Theorems 3.5, 6.1]{goodwillie3} If $F: \T \ra \Sp$ is finitary, then, for each $r$, there is a spectrum $t_r(F)$ with an action of the $r^{th}$ symmetric group $\Sigma_r$, and a natural weak equivalence
$$ D_rF(X) \simeq (t_r(F) \sm X^{\sm r})_{h\Sigma_r}.$$
\end{prop}

\begin{cor} \label{E iso cor}  If $F: \T \ra \Sp$ is finitary, and $f:X \ra Y$ is an $E_*$--isomorphism, then $P_rF(f): P_rF(X) \ra P_rF(Y)$ is also an $E_*$--isomorphism.
\end{cor}
\begin{proof} Standard spectral sequences show that any construction of the form $(C \sm X^{\sm r})_{h\Sigma_r}$ preserves $E_*$--isomorphisms in the $X$ variable. The proposition thus implies that the maps on fibers, $D_rF(f): D_rF(X) \ra D_rF(Y)$, are $E_*$--isomorphisms.  The proposition then follows by induction on $r$.
\end{proof}

\begin{rem} The corollary is false without the finitary hypothesis.  Examples can easily be constructed using homological localization functors, which are homogeneous and linear, but not, in general, finitary.
\end{rem}

\subsection{Strongly split towers of spectra}

Suppose we are given a tower of spectra under another spectrum:

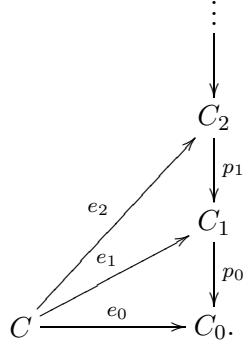
\begin{equation*}
\xymatrix{
 &  &  & \vdots \ar[d]  \\
  & &  & C_2 \ar[d]^{p_1}  \\
  &  & & C_1 \ar[d]^{p_0}   \\
 & C \ar[rr]^{e_0}  \ar[urr]^{e_1} \ar[uurr]^{e_2} & & C_0.  \\
}
\end{equation*}

We will say that the tower is {\em strongly convergent} if the connectivity of the maps $e_r$ goes to infinity as $r$ goes to infinity.

We will say that the tower is {\em strongly split} if there exists a homotopy commutative diagram

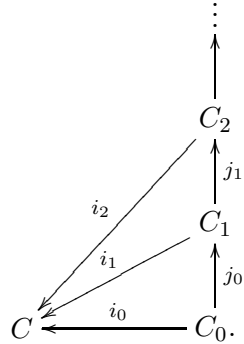
\begin{equation*}
\xymatrix{
 &  &  & \vdots   \\
  & &  & C_2 \ar[ddll]_{i_2}  \ar[u] \\
  &  & & C_1 \ar[dll]_{i_1} \ar[u]_{j_1}   \\
 & C   & & C_0 \ar[ll]_{i_0} \ar[u]_{j_0}.  \\
}
\end{equation*}

such that, for all $r$, the composite
$$ C_r \xra{i_r} C \xra{e_r} C_r$$
is an equivalence.

The following lemma is evident.

\begin{lem} \label{split tower lemma}  If a tower as above is both strongly convergent and strongly split then the induced map
$$ \hocolim_r i_r: \hocolim_r C_r \ra C$$
is an equivalence.  Thus, if $E_*$ is a homology theory, then
$$ \colim_r E_*(C_r) \ra E_*(C)$$
is an isomorphism.
\end{lem}

\begin{rem} This lemma says all that we will need to know about strongly split towers for our purposes.  However, it is illuminating to note the following.  If a tower is strongly split, one can, if needed, modify the splitting data so that $e_r \circ i_r: C_r \ra C_r$ is homotopic to the identity.  In this case, it will also true that the composite
$$ C_{r} \xra{j_r} C_{r+1} \xra{p_r} C_{r}$$
will be homotopic to the identity for all $r$.  If the tower is also strongly convergent, then there will be a wedge decomposition
$$C \simeq \bigvee_{r=1}^{\infty} \text{hofiber} \{p_r: C_r \ra C_{r-1} \}.$$
\end{rem}

\subsection{A useful proposition}

\begin{prop} \label{E iso prop} Let $F: \T \ra \Sp$ be finitary.  Suppose the Goodwillie tower of $F$ is both strongly convergent and strongly split when evaluated at a space $X$.  Then, if $f: X \ra Y$ is an $E_*$--isomorphism, then $F(f): F(X) \ra F(Y)$ is $E_*$--monic.
\end{prop}
\begin{proof}   Let $i_r: P_rF(X) \ra F(X)$ and $j_r: P_rF(X) \ra P_{r+1}F(X)$ denote the maps splitting the tower $\{P_rF(X)\}$. Suppose $f: X \ra Y$ is an $E_*$--isomorphism, and consider the diagram

\begin{equation*}
\xymatrix{
P_rF(X)  \ar[rr]^{i_r} && F(X) \ar[dd]^{F(f)}  \ar[rr]^{e_r} && P_rF(X) \ar[dd]^{P_rF(f)} \\
&&&& \\
&&  F(Y) \ar[rr]^{e_r} && P_rF(Y). } 
\end{equation*}

\vspace{.3in}

By assumption, the top composite is an equivalence, and thus an $E_*$--isomorphism. Since $f$ is an $E_*$--isomorphism, so is the right vertical map, by  \corref{E iso cor}.  We conclude that $E_*(F(f))$ is monic when restricted to the image of $E_*(i_r)$.  But $E_*(F(X))$ is the colimit over $r$ of these images, by \lemref{split tower lemma}, and thus $E_*(F(f))$ is also monic.
\end{proof}

To apply this proposition, we need criteria ensuring that a Goodwillie tower $\{P_rF(X)\}$ strongly splits.  This is the topic of our next two subsections.

\subsection{Goodwillie towers of functors with polynomial filtration}

Say that a functor $C: \T \ra \Sp$ has a {\em polynomial filtration} if it is filtered by functors $F_0C \ra F_1C \ra \dots$ such that 
$$ \hocolim_r F_rC(X) \ra C(X)$$
is an equivalence, and the homotopy cofiber functor
$$ F_rC/F_{r-1}C$$
is homogeneous of degree $r$ for all $r$.

The following lemma is well known folk knowledge.

\begin{lem} \label{poly filt lemma}In this situation, the composite
$$ F_rC(X) \ra C(X) \ra P_rC(X)$$
is an equivalence for all $r$ and $X$.  It follows that the Goodwillie tower $\{P_rC(X)\}$ will be strongly split.
\end{lem}
\begin{proof}  We have a homotopy commuative diagram with rows that are cofibration sequences
\begin{equation*}
\xymatrix{
F_rC(X) \ar[d] \ar[r] & C(X) \ar[d] \ar[r] & C/F_rC(X) \ar[d]  \\
P_rF_rC(X)  \ar[r] & P_rC(X)  \ar[r] & P_r(C/F_rC)(X). } 
\end{equation*}
As $F_rC$ has degree $r$, the left vertical map is an equivalence.  If we check that $P_r(C/F_rC)(X) \simeq *$, then the bottom left map will be an equivalence, and we will be done.  To check this we have
$$ P_r(C/F_rC)(X) \simeq \hocolim_s P_r(F_sC/F_rC)(X) \simeq *,$$
as $P_r(F_sC/F_rC)(X) \simeq *$ for $s \geq r$.
\end{proof}

\subsection{Stable natural equivalences}

Call a natural transformation $\Theta(X): C(X) \ra G(X)$ a {\em stable equivalence} if it is an equivalence for all suitably connected spaces $X$. 

\begin{lem}  If $\Theta: C \ra G$ is a stable equivalence, then 
$$ P_r\Theta(X): P_rC(X) \ra P_rG(X)$$
is an equivalence for {\em all} $X$.
\end{lem}
\begin{proof} An examination of the construction of $P_r$ shows that if $\Theta(X)$ is an equivalence for all $(d-1)$--connected spaces $X$, then so is $P_r\Theta(X): P_rC(X) \ra P_rG(X)$.  Now we can apply \cite[Corollary 3.8]{goodwillie3} which implies that, for any functor $F:\T \ra \Sp$, $P_rF$ is determined by its values on $d$--fold suspensions.
\end{proof}

\begin{cor} \label{stable nat equiv cor}Suppose $\Theta: C \ra G$ is a stable equivalence.  If, for a particular $X$, the Goodwillie tower $\{P_rC(X)\}$ is strongly split, then so also is the Goodwillie tower $\{P_rG(X)\}$.
\end{cor}

\section{Proof of \thmref{main theorem}} \label{section 3}

The Goodwillie tower of the functor from spaces to spectra sending $X$ to $\Sigma^{\infty} \Map(K,X)$ consists of a diagram of functors

\begin{equation*}
\xymatrix{
 &  &  & \vdots \ar[d]  \\
  & &  & P_3^K(X) \ar[d]  \\
  &  & & P_2^K(X) \ar[d]   \\
 & \Sigma^{\infty} \Map(K,X) \ar[rr]^{e_1}  \ar[urr]^{e_2} \ar[uurr]^{e_3} & & P_1^K(X).  \\
}
\end{equation*}
In \cite[Example 4.5]{goodwillie2}, it is shown that this tower is strongly convergent if $d(K) \leq c(X)$.  Thus \thmref{main theorem}(1) will follow from \propref{E iso prop} once we show that the tower is strongly split whenever $e(K) \leq s(X)$.  Otherwise said, we wish to show that if $n \geq e(K)$, then the tower $\{P_r^K(\Sigma^nZ)\}$ is strongly split for all spaces $Z$.  By \lemref{Pr susp lemma}, this tower agrees with the tower associated to the functor sending a space $Z$ to $\Sigma^{\infty} \Map(K,\Sigma^nZ)$.

In the terminology of the last section, the main constructions and theorems of \cite{bodigheimer} states that if $n \geq e(K)$, then there is a filtered configuration space $C(K,Z)$ such that $\Sigma^{\infty}C(K,Z)$ is a functor with a polynomial filtration, and a natural map of spaces
$$ C(K,Z) \ra \Map(K, \Sigma^nZ)$$
such that 
$$ \Sigma^{\infty} C(K,Z) \ra \Sigma^{\infty}\Map(K, \Sigma^nZ)$$
is a stable equivalence.  Then \lemref{poly filt lemma} and \corref{stable nat equiv cor} combine to say that the tower associated to $\Sigma^{\infty} \Map(K,\Sigma^nZ)$ is strongly split.

Statement (2) of \thmref{main theorem} turns out to follow easily from statement (1).  The following argument was observed by Pete Bousfield.

Suppose $g: L \ra K$ is an $E_*$--isomorphism between finite complexes. Let $X \ra X_E$ be Bousfield localization of the space $X$ with respect to $E_*$ \cite{bousfield1}. 

Consider the diagram 

\begin{equation*}
\xymatrix{
\Map(K,X) \ar[dd]^{\Map(g,X)} \ar[rr] && \Map(K,X_E)  \ar[dd]^{\Map(g,X_E)}  \\
&& \\
\Map(L,X)  \ar[rr] && \Map(L,X_E).  } 
\end{equation*}

\vspace{.1in}

As $X \ra X_E$ is an $E_*$--isomorphism, statement (1) of \thmref{main theorem} applies to say that the top map is $E_*$--monic.  The right vertical map is a homotopy equivalence as $X_E$ is $E_*$--local, and is thus an $E_*$--isomorphism.  Thus the left vertical map is $E_*$--monic. 

\begin{rem}  Though we haven't needed this here, there is an explicit model for the tower $\{P_r^K(X)\}$ for $\Sinfty \MapT(K,X)$: see \cite{arone, ahearnkuhn}.  From this model, it follows that a version of \corref{E iso cor} holds for the $K$--variable: if $E_*$ is a ring theory, and $g: L \ra K$ is an $E_*$--isomorphism, then so is $P_r^f(X): P_r^K(X) \ra P_r^L(X)$.  This leads to an alternative proof of \thmref{main theorem}(2), under the ring theory hypothesis.
\end{rem}

\appendix

\section{Computations of $e(K)$}

\begin{center} By Greg Arone and Nicholas Kuhn \end{center}

Recall that $e(K)$ is the minimal $n$ such that there exists a parallelizable $n$--dimensional manifold $M$, together with a closed subcomplex $A$ such that $K \simeq M/A$.  

In this appendix we make some observations allowing for some general estimates and explicit computations of $e(K)$.

\subsection{Upper bounds}  If $K \simeq M/A$, with $M$ parallelizable, we will say that the pair $(M,A)$ {\em represents} $K$.

\begin{lem}  $e(K\sm L) \leq e(K)+e(L)$.
\end{lem}
\begin{proof}  If $(M,A)$ represents $K$ and $(N,B)$ represents $L$, then $(M \times N, A \times N \cup M \times B)$ represents $K \sm L$.
\end{proof}

\begin{cor} \label{e(K) susp cor} $e(\Sigma^r K) \leq e(K) + r$.
\end{cor}

\begin{lem} Suppose that $L$ is a finite complex equivalent to a stably parallelizable $m$--manifold.  If $n > m$, and $K$ is obtained from $L$ by attaching $n$--cells, then $e(K) \leq m+n$.
\end{lem}
\begin{proof} Suppose $g: L \ra M$ is an equivalence, with $M$  a stably parallelizable $m$--manifold.  Let 
$$ \coprod_{i=1}^c f_i: \coprod_{i=1}^c S^{n-1} \ra L$$
be the attaching maps for constructing $K$ from $L$.  Then $\coprod_{i=1}^c S^{n-1}$ embeds in the parallelizable manifold $D^n \times M$ so that $(D^n \times M, \coprod_{i=1}^c S^{n-1})$ represents $K$.
\end{proof}

\begin{cor} \label{2 cell cor} If $K$ is the mapping cone of a map $f:S^n \ra S^m$, then $e(K) \leq m+n+1$.
\end{cor}

\begin{lem}  If $d(K) \geq 1$, then $e(K) \leq 2d(K) -1$.
\end{lem}
\begin{proof}  We can assume $K$ is a $d=d(K)$ dimensional C.W.~ complex.  Let $L$ be its $d-1$ skeleton, and let 
$$ \coprod_{i=1}^c f_i: \coprod_{i=1}^c S^{d-1} \ra L$$
denote the attaching maps of the $d$---cells of $K$.  The complex $L$ can be embedded in $\R^{2d-1}$, and we let $U$ be a regular neighborhood.  Thus $L \hra U$ is an equivalence, and $U$ is a $2d-1$ dimensional parallelizable manifold.  The composite 
$$ \coprod_{i=1}^c S^{d-1} \xra{\coprod_{i=1}^c f_i} L \hra U $$
is homotopic to an embedding, and then $(U, \coprod_{i=1}^c S^{d-1})$ will represent $K$.
\end{proof}

Slava Krushkal has told us of an unpublished result of Stallings \cite{stallings} that says that any $d$--dimensional and $c$--connected finite complex is (simple) homotopy equivalent to a subcomplex of $\mathbb R^{2d-c}$.  This implies our final upper bound.

\begin{lem} $e(K) \leq 2d(K) - c(K)$.
\end{lem}

\subsection{Lower bounds}  The obvious lower bound for $E(K)$ comes from dimension:
$$ d(K) \leq e(K).$$
Stronger lower bounds arise from the contrapositive forms of the following proposition and corollary.   

\begin{prop} \label{e(K) lower bound prop} If $n \geq e(K)$, then $\MapS(K,S^n)$ is equivalent to a suspension spectrum.
\end{prop}

Here $\MapS(K,X)$ denotes the function spectrum of stable maps between $K$ and $X$, so that  $\MapS(K,S^n)$ is the $n$--dual of $K$.

\begin{proof}  It is easy to check (see e.g. \cite{goodwillie1}) that the degree 1 approximation to $\Sinfty \MapT(K,X)$ is $\MapS(K,X)$.  Then, as in \secref{section 3}, we can conclude that, if $n \geq e(K)$ then the composite
$$ \Sinfty F_1C(K,S^0) \ra \Sinfty C(K,S^0) \ra \Sinfty \MapT(K,S^n) \ra \MapS(K,S^n)$$
is an equivalence, where $F_1C(K,S^0)$ is the first filtration of the configuration space $C(K,S^0)$.
\end{proof}

The implications of the proposition for homology are the following.

\begin{cor} \label{e(K) lower bound cor} Suppose $n \geq e(K)$. 

\noindent{\bf (1)} The reduced integral cohomology groups of $K$ satisfy
\begin{equation*}
\tilde{H}^m(K; \Z)  \text{ is } 
\begin{cases} 
0 & \text{if } m >n \\ \text{free abelian } & \text{if } m=n.
\end{cases}
\end{equation*}
\noindent{\bf (2)} For all primes $p$, $\F_p \oplus \Sigma^n\tilde{H}_{-*}(K;\F_p)$ admits the structure of an unstable algebra over the mod $p$ Steenrod algebra.
\end{cor}

For certain two cell complexes, the lower bound of the proposition matches our upper bounds.  Call a map $f: S^{m+k} \ra S^m$ {\em stably minimal} if whenever $f^{\prime}: S^{m^{\prime}+k} \ra S^{m^{\prime}}$ is a map so that $f$ and $f^{\prime}$ represent the same element in the stable homotopy group $\pi_k^S$, one has $m \leq m^{\prime}$.

\begin{thm} \label{two cell thm}  If $f: S^n \ra S^m$ is stably minimal, and $K$ is the mapping cone of $f$, then $e(\Sigma^rK) = m+n+1+r$.
\end{thm}
\begin{proof} As $S$--duality induces the identity on $\pi_*^S$, one deduces that $$\MapS(\Sigma^rK, S^{r+m+n+1}) \simeq \MapS(K,S^{m+n+1}) \simeq \Sinfty K.$$  Since $f$ is stably minimal, $\Sigma^{-i}\Sinfty K$ is not a suspension spectrum for any $i > 0$.  The proposition thus implies that $e(\Sigma^rK) \geq m+n+1+r$.  The theorem follows, as \corref{2 cell cor} and \corref{e(K) susp cor} combine to show that $e(\Sigma^rK) \leq m+n+1+r$.
\end{proof}

\subsection{Examples}  \thmref{two cell thm} applies to all the classic 2--cell complexes.  Explicitly, we have 
\begin{enumerate}
\item $e(M^{r+2}(d)) = e(\Sigma^rM^2(d)) = 3+r$,
\item $e(\Sigma^r \mathbb CP^2) = 6+r$,
\item $e(\Sigma^r \mathbb HP^2) = 12+r$,
\item $e(\Sigma^r(\text{Cayley plane})) = 24+r$, and 
\item $e(\Sigma^r(D^{2p+1} \cup_{\alpha} S^3))= 2p+4+r$ if $p$ is an odd prime and $\alpha \in \pi_{2p}(S^3)$ is an element of order $p$. 
\end{enumerate}

\end{document}